\documentclass[letterpaper, 10 pt, conference]{ieeeconf}  

\IEEEoverridecommandlockouts                              

\usepackage{amsmath}
\usepackage{amssymb}
\usepackage[svgnames]{xcolor}
\usepackage{amsfonts}
\usepackage{footnote}
\usepackage{indentfirst}
\usepackage{enumerate}
\usepackage{bbm}
\usepackage{graphicx}
\usepackage{array} 
\usepackage{adjustbox}
\usepackage{doi,color}
\usepackage{array,colortbl}
\usepackage{multicol,gensymb,multirow}
\usepackage{threeparttable,sidecap}
\usepackage{color, colortbl}
\definecolor{Gray}{gray}{0.9}

\usepackage[font=small]{caption}
\usepackage{subfigure}
\newcommand{\suchthat}{\;\ifnum\currentgrouptype=16 \middle\fi|\;}

\newcommand{\subscr}[2]{#1_{\textup{#2}}}
\newcommand{\supscr}[2]{#1^{\textup{#2}}}
\newcommand{\setdef}[2]{\{#1 \; | \; #2\}}

\newcommand{\map}[3]{#1: #2 \rightarrow #3}

\newcommand{\real}{\mathbb{R}}
\newcommand{\complex}{\mathbb{C}}

\newcommand{\imagunit}{\mathrm{i}}
  
\newcommand{\diag}{\mathrm{diag}} 

\newcommand{\fK}{\subscr{f}{K}} 
\newcommand{\fKC}{f_{\mathrm{K}\mathbb{C}}} 
\newcommand{\prj}{\mathcal{P}}
\newcommand{\bperp}{\mathrm{D}^{G}(\gamma)}
\newcommand{\scirc}{\raise1pt\hbox{$\,\scriptstyle\circ\,$}}
\newcommand{\sinc}{\mathrm{sinc}}
\newcommand\oprocendsymbol{\hbox{$\square$}}
\newcommand\oprocend{\relax\ifmmode\else\unskip\hfill\fi\oprocendsymbol}

\DeclareSymbolFont{bbold}{U}{bbold}{m}{n}
\DeclareSymbolFontAlphabet{\mathbbold}{bbold}
\newcommand{\vect}[1]{\mathbbold{#1}}

\newtheorem{theorem}{Theorem}
\newtheorem{lemma}[theorem]{Lemma}

\newtheorem{definition}[theorem]{Definition}

\newtheorem{remark}[theorem]{Remark}

\DeclareMathOperator{\Ker}{\mathrm{Ker}}
\DeclareMathOperator{\Img}{\mathrm{Img}}

\title{\LARGE \bf
Synchronization of Kuramoto Oscillators: A Power Series Approach
}

\title{\LARGE \bf
Synchronization of Coupled Oscillators:\\ The Taylor Expansion of the Inverse Kuramoto Map
}

\author{Elizabeth Y. Huang$^{1}$, Saber Jafarpour$^{1}$, Francesco Bullo$^{1}$
\thanks{*This work was supported in part by the
    U.S. Department of Energy (DOE) Solar Energy Technologies Office
    under Contract No. DE-EE0000-1583. }
\thanks{$^{1}$Department of Mechanical Engineering and the Center of Control,
  Dynamical Systems and Computation, University of California, Santa
  Barbara, 93106-5070, USA
        {\tt
    \{eyhuang, saber.jafarpour, bullo\}@engineering.ucsb.edu}}%
}

\begin{document}

\maketitle
\thispagestyle{empty}
\pagestyle{empty}

\begin{abstract}
Synchronization in the networks of coupled oscillators is a
widely studied topic in different areas. It is well-known that
synchronization occurs if the connectivity of the network
dominates heterogeneity of the oscillators. Despite extensive
study on this topic, the quest for sharp closed-form synchronization
tests is still in vain. In this paper, we present an algorithm for
finding the Taylor expansion of the inverse Kuramoto map. We show that
this Taylor series can be used to obtain a hierarchy of increasingly
accurate approximate tests with low computational complexity. These
approximate tests are then used to estimate the threshold of synchronization as well as the position of the synchronization manifold of the network.  

\end{abstract}

\section{INTRODUCTION}

\paragraph*{Problem description and motivation}

  Synchronization is a ubiquitous phenomenon which appears naturally in 
  various areas including physics, biology, and chemistry. One of the
  well-known models for studying synchronization of
  coupled oscillators is the Kuramoto model. Kuramoto model and its
  generalizations has been used for studying synchronization in various applications including 
  pacemaker cells in heart~\cite{DCM-EPM-JJ:87}, neural
  oscillators~\cite{EB-JM-PH:04}, deep brain simulation~\cite{PAT:03},
  spin glass models~\cite{GJ-JA-DB-ACCC-CPV:01}, oscillating
  neutrinos~\cite{JP:98}, chemical oscillators~\cite{IZK-YZ-JLH:02}, 
  multi-vehicle coordination~\cite{DJK-PL-KAM-TJ:08, RS-DP-NEL:07},
  synchronization of smart grids~\cite{FD-MC-FB:11v-pnas}, security
  analysis of power flow equations~\cite{AA-SS-VP:81,FW-SK:82}, optimal
  generation dispatch~\cite{JL-SHL:12}, and droop-controlled inverters
  in microgrids~\cite{MCC-DMD-RA:93, JWSP-FD-FB:12u}.

  Frequency synchronization
  is arguably one of the most important notions of synchronization for
  Kuramoto model. Frequency synchronization of Kuramoto coupled oscillators has been studied
  extensively in physics, dynamical system,
  and control communities. In the literature, many efficient
  numerical methods have been developed to study frequency 
  synchronization of the Kuramoto model. However, for the purpose of
  analysis and design of networks of coupled oscillators, these numerical
  methods have several drawbacks. First, in many
  applications (such as power networks), the network contains many
  time-varying or unknown parameters (such as power
  injections/demands). Therefore any numerical study of
  synchronization requires implementing computations for a wide range of
  parameters, which is not computationally
  tractable for large networks. Secondly, these numerical
  methods usually do not provide any guarantee for the
  performance of the synchronization such as robustness with respect to
  disturbances. Finally, these numerical methods do
  not provide any intuition about the role of different
  factors, such as network topology and dissimilarity of the oscillators, in the synchronization of Kuramoto model. 
  Therefore they are not suited for design purposes. These observations motivate
  the quest for sharp analytic conditions for synchronization of Kuramoto model~\cite{JAA-LLB-CJPV-FR-RS:05,FD-FB:13b,
  AA-ADG-JK-YM-CZ:08}. The first rigorous analytic characterization of frequency synchronization is developed
  for the complete unweighted
  graphs~\cite{DA-JAR:04,REM-SHS:05,MV-OM:08}. Further
  characterization of the frequency synchronization has been obtained
  for complete unweighted bipartite
  graphs~\cite{MV-OM:09}. Unfortunately these synchronization
  conditions require solving implicit algebraic equations and they are not
  suitable for studying the synchronization performance of the
  network. Moreover, they are only applicable to specific
  networks and cannot be extended for characterizing frequency
  synchronization of the Kuramoto model with general topologies. For the Kuramoto
  model with general topology and arbitrary weights, an ingenious
  approach based on graph theoretic ideas is proposed
  in~\cite{AJ-NM-MB:04}. Using this approach and Lyapunov analysis, a concise closed-form condition for
  synchronization of coupled oscillator is obtained~\cite{AJ-NM-MB:04}. Similar conditions have been derived in the literature using quadratic
  Lyapunov function~\cite{NC-MWS:08, FD-FB:09z} and sinusoidal
  Lyapunov function~\cite{AF-AC-WPL:10}. The essence of all these
  approaches can be explained using this intuitive idea:
  synchronization is determined by the trade-off between the connectivity of the network
  and the dissimilarity of the oscillators. Despite being elegant and
  insightful, these conditions usually provide conservative
  estimates of the synchronization threshold for large networks.  
  
  To get sharp closed-form synchronization tests for coupled
  oscillators, it is conjectured that instead of comparing the heterogeneity of
  oscillators with the well-established measures of connectivity such as spectral
  connectivity or nodal degree, one should focus on mixed measures of connectivity
  and heterogeneity of the network~\cite{FD-MC-FB:11v-pnas}. Inspired by
  characterization of frequency synchronization for acyclic networks, a suitable
  measure of connectivity and heterogeneity is proposed
  in~\cite{FD-MC-FB:11v-pnas}. Using numerical analysis for several
  random graphs and IEEE test cases, it has been shown that this condition provides a good estimate of the
  threshold of synchronization for different classes of graphs. Using cutset projections, a rigorous,
  more-conservative, and generally-applicable modified version of this
  test has been established for synchronization in~\cite{SJ-FB:16h}.

\paragraph*{Contribution} 

In this paper, we obtain a recursive scheme to compute the Taylor
expansion of the inverse Kuramoto map for the network of coupled
Kuramoto oscillators with arbitrary topology. Using results from theory of
complex analysis in several variables, we find an estimate on the domain of
convergence of this Taylor series. We show that this Taylor expansion
provides a mathematical explanation for the correct form of the
trade-off between the coupling strengths and oscillators heterogeneity
for synchronization of coupled oscillators. By truncating the Taylor series, we
obtain a family of approximate synchronization tests which can be used
to estimate the threshold of synchronization as well as find the approximate position of the synchronization manifold of the network. Using numerical
simulations, we show that these approximate tests provide a
reasonably-accurate computationally-efficient tool for estimating the synchronization of Kuramoto model.

\section{Notation}
For $n\in \mathbb{Z}_{>0}$, let $\mathbb{T}^n$ denote the $n$-torus,
let $\vect{1}_n$ (resp. $\vect{0}_n$) denote the vector
in $\real^n$ with all entries equal to $1$ (resp. $0$), and define the
vector subspace $\vect{1}_n^{\perp}=\setdef{x\in
  \real^n}{\vect{1}_n^{\top}x=0}$. Similarly, define the
vector subspace
$\vect{1}_{\complex}^{\perp}=\setdef{x\in\complex^n}{\vect{1}_n^{\top}x=0}$.
We denote $\sqrt{-1}$ by $\imagunit$. For every
$\mathbf{x}=(x_1,\ldots,x_n)^{\top}\in \real^n$, we define
$\mathrm{exp}(\imagunit \mathbf{x}) = (e^{\imagunit x_1},
\ldots, e^{\imagunit x_n})^{\top} \in \mathbb{T}^n$. A set
$R\subseteq \complex^n$ is a Reinhardt domain, if, for every
$(x_1,\ldots,x_n)^{\top}\in \real^n$ and every
$(z_1,\ldots,z_n)^{\top}\in R$, we have $(e^{\imagunit
  x_1}z_1,\ldots,e^{\imagunit x_n}z_n)^{\top}\in R$. For a complex number $z = x+ \imagunit y\in \complex$, the
norm of $z$ is denoted by $\|z\| = \sqrt{x^2+y^2}$. For a matrix $A
= \{a_{ij}\}\in \complex^{n\times m}$, the Hermitian of $A$ is denoted
by $A^H\in \complex^{m\times n}$ and the $\infty$-norm of $A$ is
defined by $\|A\|_{\infty} = \max_{i} \sum_{j=1}^{m}\|a_{ij}\|$. The
Moore\textendash{}Penrose pseudoinverse of $A$ is the unique
$A^{\dagger}\in \complex^{m\times n}$ which satisfies $AA^{\dagger}A=A$,
$A^{\dagger}AA^{\dagger}=A^{\dagger}$,
$\left(AA^{\dagger}\right)^{H}=AA^{\dagger}$, and
$\left(A^{\dagger}A\right)^{H}=A^{\dagger}A$. We define the function $\sinc:\complex\to\complex$ by 
\begin{align*}
\sinc(z) = 
\begin{cases}
\frac{\sin(z)}{z} & z\ne 0,\\ 
1 & z = 0.
\end{cases}
\end{align*}
For a complex vector $\mathbf{x}\in \complex^n$, we denote
the $\diag(\mathbf{x})$ by $[\mathbf{x}]$. For complex vectors
$\mathbf{x},\mathbf{y}\in \complex^n$ and we define their Hadamard
product by $\mathbf{x}\circ \mathbf{y} = [\mathbf{x}]\mathbf{y} =
[\mathbf{y}]\mathbf{x}$. For $n\in \mathbb{Z}_{\ge 0}$, we define the Hadamard power
of a complex vector $\mathbf{x}$ by 
\begin{align*}
(\mathbf{x})^{\circ n} = \underbrace{[\mathbf{x}]\ldots [\mathbf{x}]}_{n-1}\mathbf{x}. 
\end{align*}
We define subspaces $\Img_{\complex}(M)\subseteq \complex^m$ and $\Ker_{\complex}(M)\subseteq \complex^n$ as follows:
\begin{align*}
\Img_{\complex}(M) &=
\setdef{\mathbf{x}\in \complex^m}{\exists \ \mathbf{y}\in \complex^n \mbox{
    s.t. } \mathbf{x}=M\mathbf{y} }, \\
\Ker_{\complex}(M) &=
\setdef{\mathbf{x}\in \complex^n}{M\mathbf{x} = \vect{0}_m}.
\end{align*}
It is easy to see that the real subspaces $\Img(M) $ and $\Ker(M)$ are the
restrictions of the complex subspaces $\Img_{\complex}(M)$ and
$\Ker_{\complex}(M)$ to the real Euclidean space, respectively.  
The number of partitions of the integer $2n+1$ into $2k+1$
odd numbers is denoted by $p(k,n)$. An element in $p(k,n)$ is denoted
by $\boldsymbol{\alpha} = (\alpha_1,\ldots,\alpha_{2k+1})^{\top}$ where
$\sum_{i=1}^{2k+1}\alpha_{i} = 2n+1$. Let $G$ be a weighted undirected
connected graph with the node set $\mathcal{N}=\{1,\ldots,n\}$, the
edge set $\mathcal{E}\subseteq \mathcal{N}\times \mathcal{N}$, and the adjacency matrix $ A\in\real^{n\times n} $. We denote the incidence matrix of $G$
by $B$ and the Laplacian matrix of the graph $G$ by $L$. It is known that $L =
B\mathcal{A}B^{\top}$, where $\mathcal{A}\in \real^{m\times m}$ is the
diagonal weight matrix defined by $\mathcal{A}_{kl}=a_{ij}$ for $k=l=(i,j)$ and
$0$ otherwise. For every complex vector $\mathbf{w}\in \complex^n$, we define
$L_{\mathbf{w}} = B\mathcal{A}[\mathbf{w}]B^{\top}$.  The cutset
projection of $G$, denoted by $\prj$, is the oblique projection onto the cutset space of $G$ parallel to the cycle space of
$G$~\cite[Theorem 4]{SJ-FB:16h}. One
can show that $\prj = B^{\top}L^{\dagger}B\mathcal{A}$~\cite[Theorem
4]{SJ-FB:16h}. Additional properties of the cutset projection can be
found in~\cite{SJ-FB:16h}. Let $\subscr{T}{s}$ be a spanning tree of
the graph $G$. Then we denote the incidence matrix of $\subscr{T}{s}$
by $\subscr{B}{s}\in \real^{n\times (n-1)}$. We define the matrix
$\subscr{B^{\#}}{s}\in \real^{m\times (n-1)}$ by $\subscr{B^{\#}}{s} =
(\subscr{B^{\dagger}}{s}B)^{\top}$. Since $\subscr{T}{s}$ is a
spanning tree of $G$, we have
$\subscr{B}{s}\subscr{B^{\dagger}}{s}B = B$. This implies that
$\subscr{B^{\#}}{s}\subscr{B^{\top}}{s}= B^{\top}$. 

\section{Heterogenous Kuramoto model}

The Kuramoto model is a system of $n$ oscillators, where each
oscillator has a natural frequency $\omega_i\in \real$ and its state
is represented by a phase angle $\theta_i\in \mathbb{S}^1$. The
interconnection of these oscillators are described using a weighted
undirected connected graph $G$, with nodes
$\mathcal{N}=\{1,\ldots,n\}$, edges $\mathcal{E}\subseteq \mathcal{N}\times \mathcal{N}$, and
positive weights $a_{ij}=a_{ji}>0$ for all $ij\in\mathcal{E}$. The dynamics for the heterogeneous
Kuramoto model is given by:
\begin{equation*}
\dot{\theta}_i=\omega_i-\sum_{j=1}^{n}a_{ij}
\sin(\theta_i-\theta_j),\qquad\text{for } i\in\{1,\ldots,n\}.
\end{equation*}
In matrix language, the Kuramoto model is given by:
\begin{equation}\label{eq:kuramoto_model}
\dot{\theta}=\omega-B\mathcal{A}\sin(B^{\top}\theta),
\end{equation}
where $\theta=(\theta_1,\theta_2,\ldots,\theta_n)^{\top}\in \mathbb{T}^n$ is the phase vector,
$\omega=(\omega_1,\omega_2,\ldots,\omega_n)^{\top}\in \real^n$ is the
natural frequency vector. For every $s\in [0,2\pi)$, the clockwise
  rotation of $\theta\in \mathbb{T}^n$ by the angle $s$ is the
  function $\mathrm{rot}_s:\mathbb{T}^n\to \mathbb{T}^n$ defined by
\begin{equation*}
\mathrm{rot}_s(\theta)=(\theta_1+s,\ldots,\theta_n+s)^{\top},\qquad\text{for }
\theta\in \mathbb{T}^n.
\end{equation*}
Given $\theta\in \mathbb{T}^n$, define the equivalence class $[\theta]$ by
\begin{align*}
[\theta]=\left\{\mathrm{rot}_s(\theta)\mid s\in [0,2\pi)\right\}.
\end{align*}
The quotient space of $\mathbb{T}^n$ under the above equivalence class
is denoted by $[\mathbb{T}^n]$. If $\map{\theta}{\mathbb{R}_{\ge0}}{\mathbb{T}^n}$ is a solution for the
Kuramoto model~\eqref{eq:kuramoto_model} then, for every $s\in
[0,2\pi)$, the curve $\map{\mathrm{rot}_s(\theta)}{\real_{\ge0}}{\mathbb{T}^n}$ is also a
  solution of~\eqref{eq:kuramoto_model}. A solution
$\map{\theta}{\mathbb{R}_{\ge0}}{\mathbb{T}^n}$ of the Kuramoto
model~\eqref{eq:kuramoto_model} achieves \emph{frequency
  synchronization} if there exists a synchronous frequency $\omega_{\mathrm{syn}}\in \real$ such that 
\begin{equation*}
\lim_{t\to\infty} \dot{\theta}(t)=\omega_{\mathrm{syn}}\vect{1}_n.
\end{equation*}
By summing all the equations in~\eqref{eq:kuramoto_model}, we obtain $\omega_{\mathrm{syn}}
=\frac{1}{n}\Big(\sum_{i=1}^{n}\omega_i\Big)$. Therefore, by choosing a rotating
frame with frequency $\omega_{\mathrm{syn}}$, one can assume that $\omega\in \vect{1}^{\perp}_n$. Given natural frequency
vector $\omega\in \vect{1}^{\perp}_n$, the \emph{synchronization manifold} of the Kuramoto
model is given by 
\begin{align}\label{eq:nodal_sync}
\omega = B\mathcal{A}\sin(B^{\top}\theta).
\end{align}
In many application of the Kuramoto model, such as the power networks,
it is not only important to achieve the frequency synchronization but
it is also essential that the synchronization manifold satisfies some phase
constraints. These constraints can be understood as performance
measures for the frequency synchronization. In power network
applications, these constraints are usually imposed by the thermal and
power capacity of the lines in the network. For every $\gamma\in
[0,\frac{\pi}{2})$, we define the \emph{embedded cohesive subset}
$S^{G}(\gamma)\subseteq \mathbb{T}^n$ by
  \begin{equation*}
    S^{G}(\gamma)=\setdef{\mathrm{rot}_s(\exp(\imagunit\mathbf{x}))}{\mathbf{x}\in\bperp \text{ and }s\in [0,2\pi)},
  \end{equation*}
  where $\bperp=\{\mathbf{x}\in
  \vect{1}_n^{\perp}\mid\|B^{\top}\mathbf{x}\|_{\infty}\le
  \gamma\}$. It can be shown that $S^{G}(\gamma)$ is diffeomorphic
  with $\bperp$~\cite[Theorem 8]{SJ-FB:16h}. We refer the readers to~\cite{SJ-FB:16h} for further
  properties of the embedded cohesive subset. In this rest of this paper, we restrict
  our study to frequency synchronization of the Kuramoto model~\eqref{eq:kuramoto_model}
  in the domains $S^{G}(\gamma)$, for $\gamma\in [0,\tfrac{\pi}{2})$. We consider the state space of the Kuramoto
  model~\eqref{eq:kuramoto_model} to be $[\mathbb{T}^n]$ and we
  identify the embedded cohesive subset $S^{G}(\gamma)$ by $\bperp$.

\section{Synchronization of Kuramoto model}
We start by introducing a map which arises naturally in the study of synchronization of the Kuramoto
model~\eqref{eq:kuramoto_model}. The \emph{Kuramoto map}
$\fK:\vect{1}^{\perp}_n\to \Img(B^{\top})$ is defined by
\begin{equation}\label{def:Kuramoto-map}
\fK (\mathbf{x})=\prj\sin(B^{\top}\mathbf{x}).
\end{equation}
Note that equation~\eqref{eq:nodal_sync} can be interpreted as a nodal balance
equation. By multiplying both side of equations~\eqref{eq:nodal_sync}
by $B^{\top}L^{\dagger}$, we get 
\begin{align}\label{eq:edge_sync}
B^{\top}L^{\dagger}\omega = \prj \sin(B^{\top}\mathbf{x}) = \fK (\mathbf{x}),
\end{align}
where $B^{\top}L^{\dagger}\omega$ are the edge variables associated to
$\omega$ and equation~\eqref{eq:edge_sync} can be interpreted as the
flow balance equation. It can be shown that the frequency
synchronization in $S^{G}(\gamma)$ can be characterized by
solutions of the nodal balance equation or the flow balance
equation in $S^{G}(\gamma)$. The following theorem and its proof can be found in~\cite[Theorem 10]{SJ-FB:16h}
\begin{theorem}[\textbf{Algebraic characterization of synchronization}]\label{thm:existence_uniqueness_S}
 Consider the Kuramoto model~\eqref{eq:kuramoto_model} with $\gamma\in [0,\frac{\pi}{2})$. Then the following statements are equivalent:
\begin{enumerate}[(i)]
\item\label{p1:stable} there exists a unique locally exponentially stable synchronization
  manifold $\mathbf{x}^*$ for the Kuramoto model~\eqref{eq:kuramoto_model} in $S^{G}(\gamma)$;
\item \label{p2:nodal_balance} there exists a unique solution
  $\mathbf{x}^*$ for the node
  balance equation~\eqref{eq:nodal_sync} in $S^{G}(\gamma)$;
\item \label{p3:edge_balance} there exists a unique solution
  $\mathbf{x}^*$ for the flow balance equation~\eqref{eq:edge_sync} in $S^{G}(\gamma)$;
\end{enumerate} 
\end{theorem}\smallskip

Therefore, instead of studying frequency synchronization in
$S^{G}(\gamma)$, we focus on studying the range of the Kuramoto
map~\eqref{def:Kuramoto-map}. The following theorem shows that the
Kuramoto map is invertible on the domain  $S^{G}(\gamma)$ and
presents an iterative procedure to compute the Taylor expansion of
$\fK^{-1}(\eta)$ around $\eta=\vect{0}_m$. 

\begin{theorem}[\textbf{Inverse Kuramoto map}]\label{thm:6}
Consider the Kuramoto map defined in~\eqref{def:Kuramoto-map}. 
Let $\gamma\in [0,\tfrac{\pi}{2})$. Then the following statements hold:
\begin{enumerate}[(i)]
    \item\label{p1:invertible_real analytic} there exists a real analytic
    bijective function $\fK^{-1}:\fK(S^G(\gamma))\to S^{G}(\gamma)$ such
    that
    \begin{eqnarray*}
    \fK^{-1}\scirc \fK(\mathbf{x})&=&\mathbf{x}, \qquad\text{for } \mathbf{x}\in S^{G}(\gamma),
    \\
    \fK\scirc \fK^{-1}(\eta)&=&\eta,\qquad\text{for } \eta\in \fK(S^{G}(\gamma));
    \end{eqnarray*}
    \item\label{p2:power_series} the map $B^{\top}\fK^{-1}(\eta)$ has the
      following Taylor expansion around $\eta=\vect{0}_m$:
    \begin{align}\label{eq:4}
    \sum_{j=0}^{\infty} A_{2j+1}(\eta),
    \end{align}
    where $A_{2j+1}(\eta)$ is a homogeneous polynomial of order $2j+1$
    in $\eta$ computed iteratively by:
    \begin{align*}
    A_1(\eta) & = \eta, \\
    A_{2j+1}(\eta) & =  \prj\Bigg( \sum_{k=1}^j \frac{(-1)^{k+1}}{(2k+1)!}    
    \sum_{\boldsymbol{\alpha}\in p(k,n)} A_{\boldsymbol{\alpha}}(\eta)\Bigg),
     \end{align*}
where, for $\boldsymbol{\alpha}=(\alpha_1,\ldots,\alpha_{2k+1})\in p(k,n)$, we denote $A_{\boldsymbol{\alpha}}(\eta) = A_{\alpha_1}(\eta)\circ\ldots\circ A_{\alpha_{2k+1}}(\eta)$.
\end{enumerate}
\end{theorem}\smallskip
\begin{proof}
Regarding part~\eqref{p1:invertible_real analytic}, we first show that
$\fK$ is a local diffeomorphism on $S^{G}(\gamma)$. Note that, for
every $\mathbf{x}\in S^{G}(\gamma)$, we have $D_{\mathbf{x}} \fK =
\prj [\cos(B^{\top}\mathbf{x})] B^{\top}$. Suppose that $v\in
\vect{1}_n$ is such that $D_{\mathbf{x}} \fK(v) = \vect{0}_n$. This
implies that $[\cos(B^{\top}\mathbf{x})] B^{\top}(v)\in
\Ker(\prj) = \Ker(B\mathcal{A})$ and therefore
\begin{align*}
v^{\top}B\mathcal{A}[\cos(B^{\top}\mathbf{x})] B^{\top}(v) = 0.
\end{align*}
Since $\mathbf{x}\in S^{G}(\gamma)$, we have $\cos(x_i-x_j)>0$, for
every $(i,j)\in \mathcal{E}$. Thus, the matrix
$B\mathcal{A}[\cos(B^{\top}\mathbf{x})] B^{\top}$ is positive
semidefinite with $\vect{1}_n$ being the only eigenvector associated
to the the eigenvalue $0$. This implies that $v
=\vect{0}_n$. Therefore, for every $\mathbf{x}\in S^{G}(\gamma)$, the
map $D_{\mathbf{x}} \fK$ is an isomorphism. The Inverse Function
Theorem~\cite[Theorem 2.5.2]{RA-JEM-TSR:88} implies that $\fK$ is a
local diffeomorphism. Now we show that $\fK$ is one-to-one on the
domain $S^{G}(\gamma)$. Suppose that $\mathbf{x}_1,\mathbf{x}_2\in
S^{G}(\gamma)$ are such that $\fK(\mathbf{x}_1) =
\fK(\mathbf{x}_2)$. This means that
$B\mathcal{A}\sin(B^{\top}\mathbf{x}_1)
=B\mathcal{A}\sin(B^{\top}\mathbf{x}_2)$. As a result, we have
\begin{align*}
0 &= (\mathbf{x}_1-\mathbf{x}_2)^{\top} B\mathcal{A}
  (\sin(B^{\top}\mathbf{x}_1) - \sin(B^{\top}\mathbf{x}_2)) \\ & = (B^{\top}\mathbf{x}_1-B^{\top}\mathbf{x}_2)^{\top}\mathcal{A}
  (\sin(B^{\top}\mathbf{x}_1) - \sin(B^{\top}\mathbf{x}_2)).
\end{align*}
Since the function $y\mapsto \sin(y)$ is strictly increasing on the
interval $(-\frac{\pi}{2},\frac{\pi}{2})$, we get $B^{\top}\mathbf{x}_1=B^{\top}\mathbf{x}_2$. This implies that
$\mathbf{x}_1=\mathbf{x}_2$ and the function $\fK$ is one-to-one on
$S^{G}(\gamma)$. Note that an injective local diffeomorphism with compact domain is a 
diffeomorphism onto its image~\cite[Proposition 7.4]{JML:03}. Therefore, the map $\fK:
S^{G}(\gamma) \to \fK(S^{G}(\gamma))$ is a global diffeomorphism and part~\eqref{p1:invertible_real analytic} follows from this
result. 

Regarding part~\eqref{p2:power_series}, suppose that
$\sum_{j=0}^{\infty} A_{j}(\eta)$ is the Taylor expansion of
$B^{\top}\fK^{-1}(\eta)$ around $\eta=\vect{0}_m$. Then we have 
\begin{align*}
\eta & = \fK\scirc \fK^{-1}(\eta) = \prj \sin(B^{\top}\fK^{-1}(\eta))
  \\ & =
  \prj \sum_{k=0}^{\infty} \frac{(-1)^k}{(2k+1)!}
  \big(\sum_{i=0}^{\infty} A_{i}(\eta)\big)^{\circ (2k+1)}
\end{align*}
By equating the same order terms on both side of the above equality,
we get $A_1(\eta) =\eta$. Moreover, for every $j\in \mathbb{Z}_{\ge
  0}$, we get $A_{2j}(\eta) = \vect{0}_m$. For the odd order terms we
get
\begin{align*}
 A_{2j+1}(\eta) = \prj\Bigg( \sum_{k=1}^j \frac{(-1)^{k+1}}{(2k+1)!}  \sum_{\boldsymbol{\alpha}\in p(k,n)} A_{\boldsymbol{\alpha}}(\eta)\Bigg).
\end{align*}
This completes the proof of the theorem. 
\end{proof}

\begin{remark}
\begin{enumerate}
\item  One can compute the first five terms in the Taylor
  series~\eqref{eq:4} as follows:
\begin{align*}
A_1(\eta) &= \eta, \quad  A_3(\eta) = \prj \left(\frac{1}{3!} (\eta)^{\circ 3}\right),\\
A_5(\eta) &= \prj \left(\frac{3}{3!} A_3(\eta)\circ (\eta)^{\circ 2} -
            \frac{1}{5!} (\eta)^{\circ 5}\right).\\
A_7(\eta) & = \prj \Big(\frac{3}{3!} A_5(\eta)\circ (\eta)^{\circ 2} +
            \frac{3}{3!} (A_3(\eta))^{\circ 2}\circ \eta \\ &\quad -
            \frac{5}{5!} A_3(\eta)\circ (\eta)^{\circ 4} + \frac{1}{7!} (\eta)^{\circ 7}\Big).
\end{align*}
\item Using Theorem~\ref{thm:6}\eqref{p2:power_series}, it is easy to
  see that the Taylor expansion for $\fK^{-1}(\eta)$ is given by
  $\sum_{j=0}^{\infty} B_{2j+1}(\eta)$, where we have $ B_{2j+1}(\eta) =
  L^{\dagger}B\mathcal{A} (A_{2j+1}(\eta))$, for every $j\in \mathbb{Z}_{\ge 0}$. \oprocend
\end{enumerate}
\end{remark}\smallskip

\begin{figure}[!htb]\centering
  \includegraphics[width=.75\linewidth]{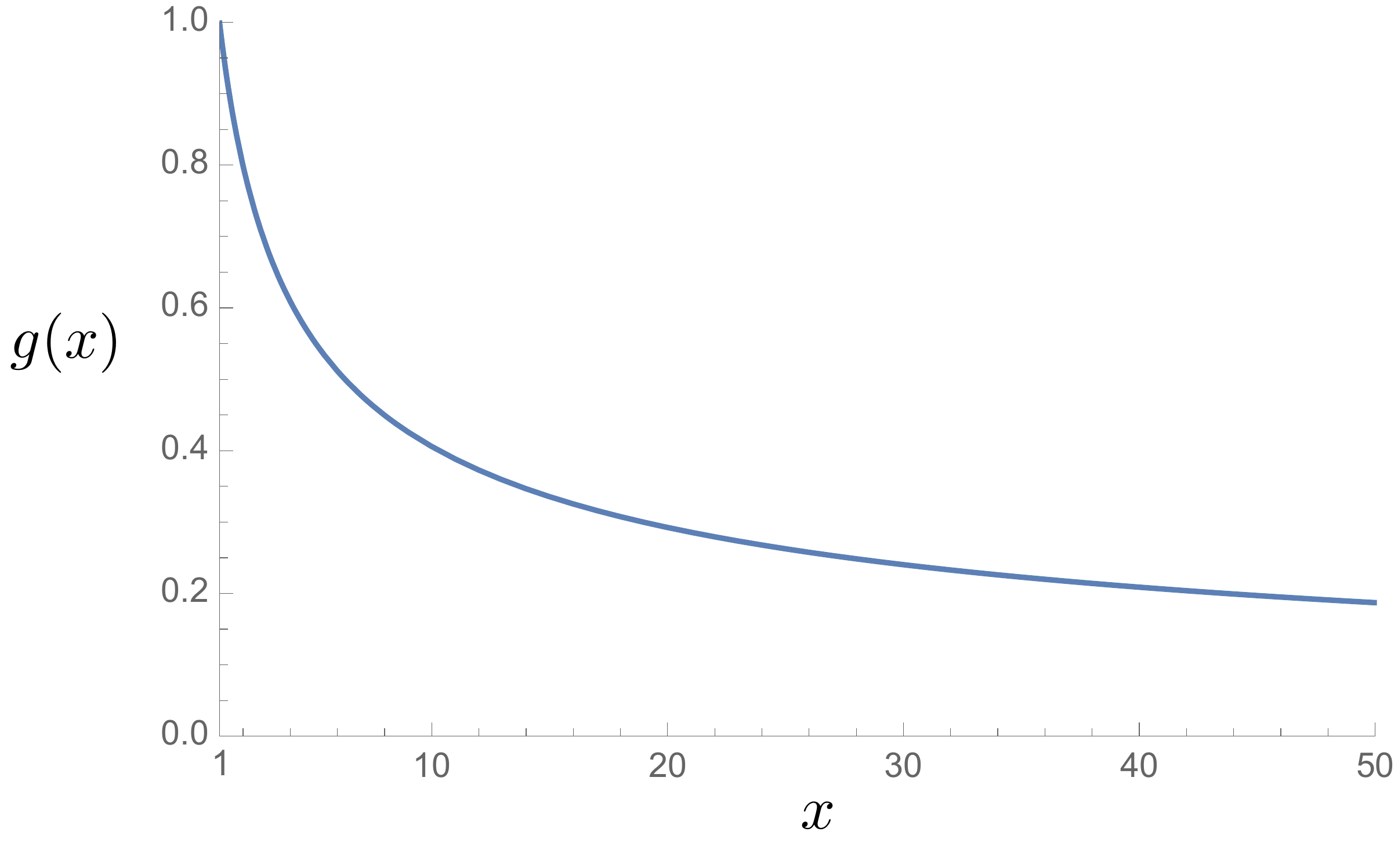}
  \caption{The graph of the monotonically-decreasing function $g$}
  \label{fig:function-g}
\end{figure}

Before we state the main result of this paper, it is convenient to introduce the smooth function
$\map{g}{[1,\infty)}{\real}$ defined by
  \begin{multline*}
    g(x)= \frac{y(x)+\sin(y(x))}{2} \\ - x \frac{y(x)-\sin(y(x))}{2}
    \,\Big|_{y(x) = \arccos(\frac{x-1}{x+1})}.
  \end{multline*}
One can verify that $g(1)=1$, $g$ is monotonically decreasing,
  and $\lim_{x\to\infty} g(x)=0$; the graph of $g$ is shown in
  Figure~\ref{fig:function-g}.
Suppose that $G$ is a weighted undirected connected graph and
$\subscr{T}{s}$ is a spanning tree of $G$. Then using function $g$, one can define two subsets of $\Img(B^{\top})$ as follows:
\begin{align*}
\Omega &= \setdef{\eta\in\Img(B^{\top})}{\|\eta\|_{\infty}\le g(\|\prj\|_{\infty})},\\
\supscr{\Omega}{s} &= \setdef{\xi\in\Img(\subscr{B^{\top}}{s})}{\|\xi\|_{\infty}\le \big\|\subscr{B^{\#}}{s}\big\|^{-1}_{\infty}g(\|\prj\|_{\infty})},
\end{align*}
It is clear that we have $\subscr{B^{\#}}{s}\supscr{\Omega}{s}\subseteq \Omega$.
In order to study the convergence of the Taylor series~\eqref{eq:4}, we need to use
various tools form theory of complex analysis in several variables. These tools are
developed in the Appendix. In the next theorem, using the results in
the Appendix, we prove the convergence of the Taylor expansion~\eqref{eq:4} on a suitable domain.

\begin{theorem}[\textbf{Synchronization of Kuramoto model}]\label{thm:sync}
Consider the Kuramoto model~\eqref{eq:kuramoto_model}, with graph $G$
and cutset projection $\prj$. Suppose that $\subscr{T}{s}$ is a
spanning tree of $G$ and define $\gamma^*\in [0,\tfrac{\pi}{2})$ by
\begin{align*}
\gamma^* = \arccos\left(\frac{\|\prj\|_{\infty}-1}{\|\prj\|_{\infty}+1}\right).
\end{align*}
Then the following statements holds: 
\begin{enumerate}[(i)]
\item\label{p1:solution} For every $\omega\in \vect{1}^{\perp}_n$ such that
  $B^{\top}L^{\dagger}\omega\in \Omega$, there exists a unique locally
  exponentially stable synchronization manifold $\mathbf{x}^*$ for the Kuramoto model~\eqref{eq:kuramoto_model} in $S^{G}(\gamma^*)$; 
\item\label{p2:series_converges} For every $\omega\in \vect{1}^{\perp}_n$ such that
  $\subscr{B^{\top}}{s}L^{\dagger}\omega\in \supscr{\Omega}{s}$, the Taylor expansion~\eqref{eq:4}:
\begin{align*}
\sum_{j=0}^{\infty} A_{2j+1}(B^{\top}L^{\dagger}\omega)
\end{align*}
 converges absolutely to $B^{\top}\mathbf{x}^*$.
\end{enumerate}
\end{theorem}\smallskip
\begin{proof}
Regarding part~\eqref{p1:solution}, by
Theorem~\ref{thm:range_kuramoto}\eqref{p2:real}, there exists
$\mathbf{x}^*\in S^{G}(\gamma^*)$ such that
$B^{\top}L^{\dagger}\omega=\fK(\mathbf{x}^*)$. The result follows from
Theorem~\ref{thm:existence_uniqueness_S}\eqref{p3:edge_balance}. 

Regarding part~\eqref{p2:series_converges}, we use the notations and
results in the Appendix. First note that by
Theorem~\ref{thm:Kuramoto_complex}\eqref{p3:invertible_analytic}, the
inverse of the complex Kuramoto map $\fKC^{-1}$ is a well-defined holomorphic diffeomorphism on $\Omega_{\complex}$. It is easy to see that
$\subscr{B^{\#}}{s}\supscr{\Omega_{\complex}}{s}\subseteq 
\Omega_{\complex}$. Therefore, one can define the holomorphic mapping
$B^{\top}\fKC^{-1}\scirc \subscr{B^{\#}}{s} :
\supscr{\Omega_{\complex}}{s}\to \Img_{\complex}(B^{\top})$ by 
\begin{align*}
B^{\top}\fKC^{-1}\scirc \subscr{B^{\#}}{s} (\mathbf{x}) =
  B^{\top}\fKC^{-1}(\subscr{B^{\#}}{s} (\mathbf{x})). 
\end{align*}
It is easy to see that $\sum_{j=0}^{\infty} A_{2j+1}(\subscr{B^{\#}}{s} (\mathbf{x}))$ is the Taylor expansion of the mapping $B^{\top}\fKC^{-1}\scirc \subscr{B^{\#}}{s} (\mathbf{x})$
around $\mathbf{x}=\vect{0}_{n-1}$. Moreover, the set
$\supscr{\Omega_{\complex}}{s}$ is a Reinhardt domain in $\complex^{n-1}$. Therefore, by~\cite[Theorem 2.4.5]{LH:90},
the Taylor series $\sum_{j=0}^{\infty} A_{2j+1}(\subscr{B^{\#}}{s} (\mathbf{x}))$ converges absolutely
to $B^{\top}\fKC^{-1}\scirc \subscr{B^{\#}}{s}(\mathbf{x})$, for every
$\mathbf{x}\in \supscr{\Omega_{\complex}}{s}$. Since $\subscr{B^{\top}}{s}L^{\dagger}\omega\in \supscr{\Omega_{\complex}}{s}$, the Taylor series
\begin{align*}
\sum_{j=0}^{\infty} A_{2j+1}(\subscr{B^{\#}}{s} ( \subscr{B^{\top}}{s}L^{\dagger}\omega))
  = \sum_{j=0}^{\infty} A_{2j+1}(B^{\top}L^{\dagger}\omega)
\end{align*}
converges absolutely to $B^{\top}\fKC^{-1}\scirc \subscr{B^{\#}}{s}(\subscr{B^{\top}}{s}L^{\dagger}\omega) = B^{\top}\fKC^{-1}(B^{\top}L^{\dagger}\omega) =
B^{\top}\mathbf{x}^*$. 
\end{proof}
Theorem~\ref{thm:sync}\eqref{p2:series_converges} provides an estimate on the domain of convergence of the power
series~\eqref{eq:4}. For many large networks, this
estimate is overly conservative and the domain of convergence is usually much larger than $\supscr{\Omega}{s}$. However, one can still use the formal
power series~\eqref{eq:4} as a tool for studying approximate
synchronization of the Kuramoto model as well as the approximate position of the
synchronization manifold. 

\begin{definition}[\textbf{Approximate synchronization}]
For $n\in \mathbb{Z}_{\ge 0}$, we define the $(2n+1)$th approximate synchronization
manifold $S_{2n+1}\in \vect{1}_n^{\perp}$ of the Kuramoto model~\eqref{eq:kuramoto_model} by 
\begin{align}\label{sol:approx}
S_{2n+1} = L^{\dagger}B\mathcal{A}\Big(\sum_{j=1}^{n}A_{2j+1}(B^{\top}L^{\dagger}\omega)\Big),
\end{align}
and the $(2n+1)$th approximate synchronization test $T_{2n+1}$
for synchronization of the Kuramoto model~\eqref{eq:kuramoto_model} in $S^{G}(\gamma)$ by
\begin{align}\label{test:approx}
\left\|\sum_{j=1}^{n}
  A_{2j+1}(B^{\top}L^{\dagger}\omega)\right\|_{\infty}\le \gamma. 
\end{align}
\end{definition}

\section{Simulations}

In this section, we study the accuracy and computational time of the family of approximate
tests~\eqref{sol:approx} and~\eqref{test:approx} on several IEEE test cases. An IEEE test case is described by a connected
graph $G$ and a nodal admittance matrix $Y\in \mathbb{C}^{n\times
  n}$. The set of nodes of $G$ is partitioned into a set of load buses
$\mathcal{N}_1$ and a set of generator buses $\mathcal{N}_2$. The
voltage at the node $j\in \mathcal{N}_1\cup\mathcal{N}_2$ is denoted
by $V_j$, where $V_j= |V_j|e^{\imagunit\theta_j}$ and the power demand
(resp. power injection) at node $j\in \mathcal{N}_1$ (resp. $j\in \mathcal{N}_2$)
is denoted by $P_j$. By ignoring the resistances in the network and
the power balance equations for reactive power, the
synchronization manifold $[\theta]$ of the network satisfies the
following Kuramoto model~\cite{FD-MC-FB:11v-pnas}:
\begin{align}\label{eq:power-network}
P_j - \sum\nolimits_{l\in \mathcal{N}_1\cup \mathcal{N}_2} a_{jl}\sin(\theta_j-\theta_l)=0,\quad\forall
  j\in \mathcal{N}_1\cup \mathcal{N}_2. 
\end{align}
where $a_{jl}=a_{lj}=|V_j||V_l|\mathrm{Im}(Y_{jl})$, for connected nodes $j$
and $l$. We consider effective power injections (demands) to be a scalar multiplication of
nominal power injections (demands), i.e., given nominal power profile
$P^{\mathrm{nom}}$ we set $P_j=K P^{\mathrm{nom}}_j$, for some $K\in
\real_{>0}$ and for every $j\in \mathcal{N}_1\cup\mathcal{N}_2$. The
voltage magnitudes at the generator buses are pre-determined and the
voltage magnitudes at load buses are computed by solving the reactive
power balance equations using the optimal power flow solver provided
by MATPOWER~\cite{RDZ-CEMS-RJT:11}. For six IEEE test cases given in
Table~\ref{tab:IEEE-test-cases}, we numerically compute the accuracy of
the approximate synchronization test $T_{5}$ in the domains
$S^{G}({\pi}/{4})$ and $S^{G}({\pi}/{2})$.
\begin{center}
\begin{table}[htb] 
\resizebox{0.99\linewidth}{!}{
\begin{tabular}{|l|c|c|c|c|}
\hline
\multirow{3}{*}{Test Case}  
			&\multicolumn{2}{c|}{$T_5$ for
                          $\gamma=\frac{\pi}{4}$} & \multicolumn{2}{c|}{$T_5$ for $\gamma=\frac{\pi}{2}$}\\
                         \cline{2-5}
  & Accuracy\textsuperscript{*}  &Computation & Accuracy &Computation \\

  & &time \textsuperscript{\dag} & & time  \\
\hline
\rowcolor{Gray}
Chow 9  & 99.57~$\%$& 0.51~$\%$& 90.05~$\%$&  0.48~$\%$ \\

IEEE 14  & 99.80~$\%$& 0.26~$\%$& 94.96~$\%$& 0.43~$\%$\\

\rowcolor{Gray}
IEEE 30 &99.80~$\%$&0.64~$\%$&93.84~$\%$& 0.82~$\%$\\

IEEE 57 &99.80~$\%$&0.83~$\%$&94.77~$\%$&1.67~$\%$\\

\rowcolor{Gray}
IEEE 118 & 99.69~$\%$&2.47~$\%$&93.60~$\%$&4.39~$\%$\\

IEEE 300 &
           99.42~$\%$&10.60~$\%$&84.46~$\%$&17.55~$\%$\\
\hline
\end{tabular}}
\begin{tablenotes}\footnotesize
                          \item \textsuperscript{*} Accuracy 
                            is compared with the solution
                            of \emph{fsolve}.
                          \item \textsuperscript{\dag} Computation
                            time is compared with the computation time
                            of \emph{fsolve}.
\end{tablenotes}
\caption{Comparison of accuracy and computation time of approximate synchronization
  test $T_5$ on six IEEE test cases.}\label{tab:IEEE-test-cases}
\end{table}
\end{center}
\begin{remark}
\begin{enumerate}
\item By increasing the value of $\gamma$, the accuracy of the
  approximate tests decreases; 
\item By increasing the number of nodes in the IEEE test case, the
  computation time of the approximate test increases. This is mainly because
  computing the approximate test requires multiplication by the matrix
  $\prj\in \real^{m\times m}$ which is a dense matrix with $n-1\le m \le
  {n \choose 2}$.  \oprocend
\end{enumerate}
\end{remark}\smallskip

\begin{figure}[htp]
	\centering
	\includegraphics[width=0.48\textwidth]{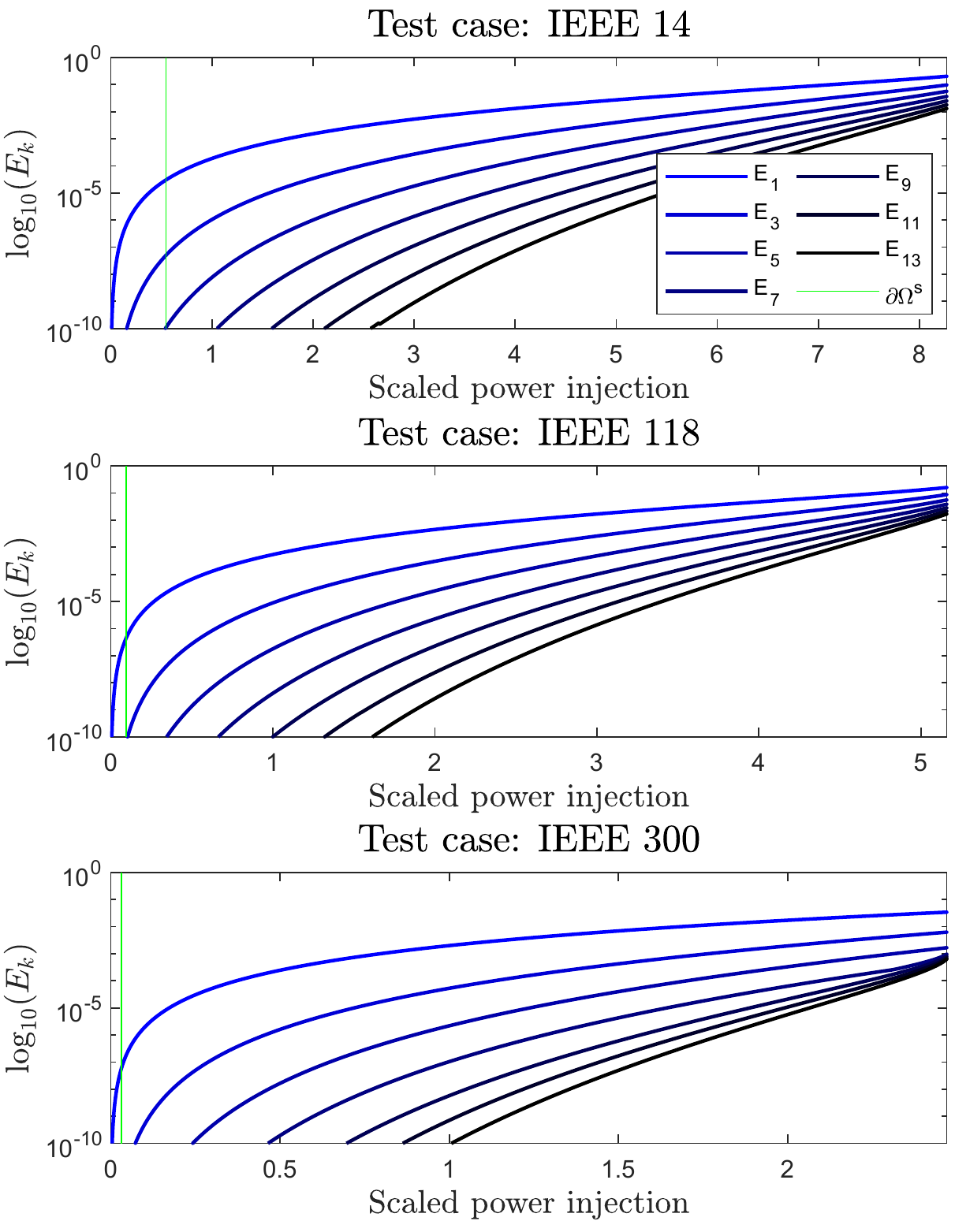}
	\caption{Comparison of accuracy of the approximate synchronization
		manifold $S_{2n+1}$ for select IEEE test cases. The green line in these
		graphs shows the boundary of the region $\supscr{\Omega}{s}$.}\label{fig:IEEE-S-test}
\end{figure}

For a select set of the IEEE test cases, we compute the accuracy of
equations~\eqref{sol:approx} in approximating the position of the
synchronization manifold of the system. For every $n\in
\mathbb{Z}_{\ge 0}$, we define the error of the approximate manifold
$S_{2n+1}$ by 
\begin{align*}
\mathrm{E}_{2n+1} = \left\|S_{2n+1} - \theta^*\right\|_{\infty}
\end{align*}
where $\theta^*\in \vect{1}^{\perp}_n$ is the solution for the power flow equations~\eqref{eq:power-network}
using \emph{fsolve} in MATLAB. The comparison of the error for
different approximate manifolds is shown in
Figure~\eqref{fig:IEEE-S-test}. 

\section{Conclusion}

We presented an algorithm to compute the Taylor expansion of the
inverse Kuramoto map. We showed that this Taylor series can be used to
study the frequency synchronization as well as to find the position of the synchronization manifold of the
Kuramoto model of coupled oscillators. Using techniques from complex
analysis in several variable, we found an estimate on the domain of
convergence of this Taylor series. Numerical simulations show
that this power series gives approximations of the
synchronization of Kuramoto model with reasonable accuracy and low
computational cost, way beyond our estimate on its domain of convergence.

\section*{Appendix}\label{app:scv}
In this appendix, we generalize the results in~\cite{SJ-FB:16h} to the complexified
Kuramoto model. We start by defining the domains
\begin{align*}
S_{\complex}^{G}(\gamma) &=\setdef{\mathbf{x}\in
                           \vect{1}_{\complex}^{\perp}}{\|B^{\top}\mathbf{x}\|_{\infty}
                           \le \gamma},\\
\Omega_{\complex} &=\setdef{\eta\in
                           \Img_{\complex}(B^{\top})}{\|\eta\|_{\infty}
                           \le g(\|\prj\|_{\infty})},\\
\supscr{\Omega_{\complex}}{s}&=\setdef{\xi\in
                           \Img_{\complex}(\subscr{B^{\top}}{s})}{\|\xi\|_{\infty}
                           \le \big\|\subscr{B^{\#}}{s}\big\|^{-1}_{\infty}g(\|\prj\|_{\infty})},
\end{align*}
and complex Kuramoto map $\fKC: \vect{1}_{\complex}^{\perp}\to \Img_{\complex}(B^{\top})$ by 
\begin{align}\label{eq:complex_kuramoto}
\fKC(\mathbf{z}) = \prj\sin(B^{\top}\mathbf{z}).
\end{align}
We first state a useful lemma. 
\begin{lemma}\label{lem:rank}
Let $\mathbf{z}\in S^{G}_{\complex}(\frac{\pi}{2})$. Then we have 
\begin{align*}
\mathrm{rank}(L_{\cos(\mathbf{z})})=\mathrm{rank}(L_{\sinc(\mathbf{z})})=n-1.
\end{align*}
\end{lemma}\smallskip
\begin{proof}
Let $\mathbf{w}\in \complex^n$. Then it is clear that $L_{\mathbf{w}}\vect{1}_n = \vect{0}_n$. Suppose that, for some $v\in
\vect{1}_{\mathbb{C}}^{\perp}$, we have $L_{\mathbf{w}}v=\mathbf{0}$. This
implies that
$B\mathcal{A}[\mathbf{w}]B^{\top}v=\vect{0}_n$. Therefore, we have $v^HB\mathcal{A}[\mathbf{w}]B^{\top}v=\vect{0}_n$.
Since $B$ and $\mathcal{A}$ are real matrices, we have 
\begin{align*}
\mathrm{Re}(v^HB\mathcal{A}[\mathbf{w}]B^{\top}v) = v^HB\mathcal{A}[\mathrm{Re}(\mathbf{w})]B^{\top}v=\vect{0}_n.
\end{align*}
Suppose that $\mathbf{z}=\mathbf{x}+\imagunit \mathbf{y}$. If we set $\mathbf{w} =
\cos(\mathbf{z})$, then $\mathrm{Re}(\mathbf{w}) =
[\cos(\mathbf{x})]\cosh(\mathbf{y})$. Therefore, for every $z\in S^{G}_{\mathbb{C}}(\frac{\pi}{2})$, the matrix
$B\mathcal{A}[\cos(B^{\top}\mathbf{x})][\cosh(B^{\top}\mathbf{y})]B^{\top}$
is positive semidefinite and its only eigenvector associated to the
eigenvalue $0$ is $\vect{1}_{n}$. Thus $v\in
\mathrm{span}_{\mathbb{C}}(\vect{1}_n)$ and we have $\mathrm{rank}(L_{\cos(\mathbf{z})})=n-1$. If we set $\mathbf{w} = \sinc(\mathbf{z})$, then we know that, for
every $i\in \{1,\ldots,m\}$, we have
\begin{multline*}
\mathrm{Re}(\sinc(z_i)) \\ = \frac{x_i\sin(x_i)\cosh(y_i) +
  y_i\cos(x_i)\sinh(y_i)}{\sqrt{x_i^2+y_i^2}}.
\end{multline*}
Therefore, for every $\mathbf{z}\in S^{G}_{\complex}(\frac{\pi}{2})$,
the matrix $B\mathcal{A}[\mathrm{Re}(\sinc(\mathbf{z}))]B^{\top}$ is positive
semidefinite and its only eigenvector associated to the
eigenvalue $0$ is $\vect{1}_{n}$. Thus $v\in
\mathrm{span}_{\mathbb{C}}(\vect{1}_n)$ and we have $\mathrm{rank}(L_{\sinc(\mathbf{z})})=n-1$.
\end{proof}

For every $\mathbf{y}\in \Img_{\complex}(B^{\top})$, we define the map $Q_{\mathbf{y}}:\Img_{\complex}(B^{\top})\to \Img_{\complex}(B^{\top})$ by 
\begin{align*}
Q_{\mathbf{y}} \mathbf{x} =\prj
  [\sinc(\mathbf{y})] \mathbf{x}, \qquad\mbox{ for all }
  \mathbf{x}\in \Img_{\complex}(B^{\top}).
\end{align*}
Then we have the following lemma.
\begin{lemma}\label{lem:Q}
Suppose that $\mathbf{z}\in \Img_{\complex}(B^{\top})$ is such that
$\|\mathbf{z}\|_{\infty}\le \frac{\pi}{2}$, then the
following statements hold: 
\begin{enumerate}[(i)]
\item\label{p1:invertible} the map $Q_{\mathbf{z}}$ is invertible;
\item\label{p2:inverse} we have $Q^{-1}_{\mathbf{z}} =
  B^{\top}L^{\dagger}_{\sinc(\mathbf{z})}B\mathcal{A}$;
\item\label{p3:cont} the map $\mathbf{z}\to Q^{-1}_{\mathbf{z}}$ is continuous on
  $\|\mathbf{z}\|_{\infty}\le \frac{\pi}{2}$. 
\end{enumerate}
\end{lemma}\smallskip
\begin{proof}
Regarding part~\eqref{p1:invertible}, to show that $Q_{\mathbf{z}}$ is invertible, it suffices to show that
it is injective and surjective. We first show that $Q_{\mathbf{z}}$ is
injective. Suppose that there exists $v\in \Img_{\complex}(B^{\top})$
such that $Q_{\mathbf{z}}v=\vect{0}_m$. Let $w\in\vect{1}_{\complex}^{\perp}$ be such that $B^{\top}w=v$. Using the fact that $\Ker_{\complex}(\prj)=\Ker_{\complex}(B\mathcal{A})$, we get
\begin{align*}
B\mathcal{A}[\sinc(\mathbf{z})]B^{\top}w=\vect{0}_n.
\end{align*}
Using Lemma~\ref{lem:rank}, we get $w=\vect{0}_n$ and
$v=\vect{0}_m$. The surjective of $Q_{\mathbf{z}}$ follows from
Rank\textendash{}nullity Theorem. 

Regarding part~\eqref{p2:inverse}, for every $v\in \Img_{\complex}(B^{\top})$, we have 
\begin{align*}
Q_{\mathbf{z}}B^{\top}L^{\dagger}_{\sinc(\mathbf{z})}B\mathcal{A}v  = \prj
  [\sinc(\mathbf{z})]
  B^{\top}L^{\dagger}_{\sinc(\mathbf{z})}B\mathcal{A}v 
\end{align*}
Using Lemma~\ref{lem:rank}, $\mathrm{rank}(L_{\sinc(\mathbf{z})}) = n-1$. This implies that
$L^{\dagger}_{\sinc(\mathbf{z})}L_{\sinc(\mathbf{z})} =
L_{\sinc(\mathbf{z})}L^{\dagger}_{\sinc(\mathbf{z})} = I_n-\frac{1}{n}\vect{1}_n\vect{1}_n^{\top}$. Therefore
\begin{align*}
Q_{\mathbf{z}}B^{\top}L^{\dagger}_{\sinc(\mathbf{z})}B\mathcal{A}v  =
  B^{\top}L^{\dagger}B\mathcal{A}v = v. 
\end{align*}
Using part~\eqref{p1:invertible}, one can deduce that
$Q_{\mathbf{z}}^{-1} =
B^{\top}L^{\dagger}_{\sinc(\mathbf{z})}B\mathcal{A}$.

Regarding part~\eqref{p3:cont}, the result follows from
part~\eqref{p2:inverse} and~\cite[Theorem 4.2]{VR:97}.
\end{proof}

Now we can prove the main result of this appendix. 

\begin{theorem}[\textbf{Range of the complex Kuramoto map}]\label{thm:range_kuramoto}
Define $\gamma^*\in [0,\frac{\pi}{2})$ by 
\begin{align*}
\gamma^* = \arccos\left(\frac{\|\prj\|_{\infty}-1}{\|\prj\|_{\infty}+1}\right).
\end{align*}
and suppose that $\eta\in \Img_{\complex}(B^{\top})$ is such that
\begin{align*}
\|\eta\|_{\infty} \le g(\|\prj\|_{\infty}).
\end{align*}
Then the following statements hold: 
\begin{enumerate}[(i)]
\item\label{p1:complex} There exists a unique $\mathbf{z}^*\in S^{G}_{\complex}(\gamma^*)$ such that $\eta = \fKC(\mathbf{z}^*)$; 
\item\label{p2:real} If $\eta\in \real^m$, then there exists a unique $\mathbf{x}^*\in S^{G}(\gamma^*)$ such that $\eta = \fKC(\mathbf{x}^*)$; 
\end{enumerate}
\end{theorem}\smallskip
\begin{proof}
Regarding part~\eqref{p1:complex}, define the set $\Gamma = \{\mathbf{x}\in \Img_{\complex}(B^{\top})\mid
\|\mathbf{x}\|_{\infty}\le \gamma^*\}$. Then one can consider the following fixed-point problem on $\Gamma$: 
\begin{align}\label{eq:fixed-point}
\mathbf{y}= Q^{-1}_{\mathbf{z}}(\eta) := h_{\eta}(\mathbf{z}). 
\end{align}
By Lemma~\ref{lem:Q}, the map $\mathbf{z}\to h_{\eta}(\mathbf{z})$ is
well-defined and continuous on $\Gamma$. We show that
$h_{\eta}(\Gamma)\subseteq \Gamma$. Let $\mathbf{z}\in \Gamma$. Then we have 
\begin{align*}
\|h_{\eta}(\mathbf{z})\|_{\infty} =
  \|Q^{-1}_{\mathbf{z}}(\eta)\|_{\infty} \le \|Q^{-1}_{\mathbf{z}}\|_{\infty}\|\eta\|_{\infty}. 
\end{align*}
By~\cite[Lemma 24]{SJ-FB:16h}, we have 
\begin{align*}
  \|Q^{-1}_{\mathbf{z}}\|_{\infty} = \left(\min_{\mathbf{z}\in \Gamma}
  \min_{\substack{\mathbf{x}\in \Img_{\complex}(B^{\top}), \\ \|\mathbf{x}\|_{\infty}=1}} \|Q_{\mathbf{z}}\mathbf{x}\|_{\infty}\right)^{-1}.
\end{align*}
Let $\mathbf{z}\in \Gamma$, $\mathbf{x}\in
\Img_{\complex}(B^{\top})$ be such that $\|\mathbf{x}\|_{\infty}=1$, and
$\mathbf{w} = \sinc(\gamma^*)\vect{1}_n$. Then we have  
\begin{align*}
\|Q_{\mathbf{y}}\mathbf{x}\|_{\infty} &=
  \|\prj[\sinc(\mathbf{z})]\mathbf{x}\|_{\infty} \\ & \ge
  \|\prj[\mathbf{w}]\mathbf{x}\|_{\infty} - \|\prj[\mathbf{w}-\sinc(\mathbf{z})]\mathbf{x}\|_{\infty}.
\end{align*}
The first term on the right hand side is given by:
\begin{align*}
  \|\prj[\mathbf{w}]\mathbf{x}\|_{\infty}=\frac{1+\sinc(\gamma^*)}{2}
\end{align*}
The second term on the right hand side can be bounded as:
\begin{align*}
\|\prj[\mathbf{w}-\sinc(\mathbf{z})]\mathbf{x}\|_{\infty}\le
  \|\prj\|_{\infty}\|[\mathbf{w}-\sinc(\mathbf{z})]\|_{\infty} 
\end{align*}
Note that, if $\mathbf{z} = \mathbf{x} + \imagunit \mathbf{y}$, then
we have 
\begin{multline*}
\|[\mathbf{w}-\sinc(\mathbf{z})]\|_{\infty} \\=
\max_{i} \left|\frac{1+\sinc(\gamma^*)}{2} -
  \frac{\sin(x_i)^2+\sinh(y_i)^2}{\sqrt{x^2_i+y^2_i}}\right|.
\end{multline*}
Therefore,
\begin{align*}
\|\prj[\mathbf{w}-\sinc(\mathbf{z})]\mathbf{x}\|_{\infty}\le \|\prj\|_{\infty}\frac{1-\sinc(\gamma^*)}{2}
\end{align*}
This implies that
\begin{multline*}
\min_{\mathbf{z}\in \Gamma} \min_{\substack{\mathbf{x}\in
  \Img_{\complex}(B^{\top}), \\ \|\mathbf{x}\|_{\infty}=1}}
  \|Q_{\mathbf{z}}\mathbf{x}\|_{\infty} \\ \ge
  \frac{1+\sinc(\gamma^*)}{2}-\|\prj\|_{\infty}\frac{1-\sinc(\gamma^*)}{2}\\
  = \frac{g(\|\prj\|_{\infty})}{\gamma^* },
\end{multline*}
where the last equality is by definition of function $g$. Thus
\begin{align*}
\|h_{\eta}(\mathbf{z})\|_{\infty} \le
  \|Q^{-1}_{\mathbf{z}}\|_{\infty}\|\eta\|_{\infty} \le \gamma^*.
\end{align*}
Therefore, by the Brouwer Fixed-point Theorem, the map $h_{\eta}$ has a
fixed point in $\Gamma$, i.e., there exists $\mathbf{y}^*\in \Gamma$ such that 
\begin{align*}
\mathbf{y}^* = h_{\eta}(\mathbf{y}^*) = Q^{-1}_{\mathbf{y}^*}\eta.
\end{align*}
This means that we have 
\begin{align*}
\eta = Q_{\mathbf{y}^*}\mathbf{y}^* = \prj \sin(\mathbf{y}^*). 
\end{align*}
Since $\mathbf{y}^*\in \Gamma$, there exists $\mathbf{z}^*\in S^{G}_{\complex}(\gamma^*)$ such that $\mathbf{y}^* =
B^{\top}\mathbf{z}^*$. This completes the proof of part~\eqref{p1:complex}.

Regarding part~\eqref{p2:real}, the proof is given in~\cite[Theorem 16 (ii)]{SJ-FB:16h}.
\end{proof}

Using Theorem~\ref{thm:range_kuramoto}, one can show that the complex Kuramoto map $\fKC$ is a
diffeomorphism on $\Omega_{\complex}$. 

\begin{theorem}[\textbf{Inverse complex Kuramoto map}]\label{thm:Kuramoto_complex}
Consider the complex Kuramoto map defined in~\eqref{eq:complex_kuramoto}. 
Let $\gamma\in [0,\tfrac{\pi}{2})$. Then the following statements hold:
\begin{enumerate}[(i)]
    \item\label{p1:local-diff} $D_{\mathbf{z}}\fKC$ is invertible, for every $\mathbf{z}\in S_{\complex}^{G}(\gamma)$;
    \item\label{p3:invertible_analytic} There exists a holomorphic
      bijective function $\fKC^{-1}:\Omega_{\complex}\to \fKC^{-1}(\Omega_{\complex})$ such
    that
    \begin{eqnarray*}
    \fKC^{-1}\scirc \fKC(\mathbf{x})&=&\mathbf{x}, \qquad\text{for } \mathbf{x}\in
                                    \fKC^{-1}(\Omega_{\complex}),
    \\
    \fKC\scirc \fKC^{-1}(\eta)&=&\eta,\qquad\text{for } \eta\in \Omega_{\complex};
    \end{eqnarray*}
\end{enumerate}
\end{theorem}
\smallskip
\begin{proof}
Regarding part~\eqref{p1:local-diff}, pick $\mathbf{z}\in
S^{G}_{\complex}(\gamma)$. First note that, we have
\begin{align*}
D_{\mathbf{z}}\fKC = \prj[\cos(B^{\top}\mathbf{z})]B^{\top}.
\end{align*}
Since $\Img_{\complex}(\prj)=\Img_{\complex}(B^{\top})$, we know that
$\mathrm{rank}(D_{\mathbf{z}}\fKC) =
\mathrm{rank}(L_{\cos(\mathbf{z})})$. The result then follows from
Lemma~\ref{lem:rank}. 

Regarding part~\eqref{p3:invertible_analytic}, the
restriction of the function
$\fKC$ to the open set $\fKC^{-1}(\Omega_{\complex})$ is a holomorphic map
$\fKC\mid_{\fKC^{-1}(\Omega_{\complex})} : \fKC^{-1}(\Omega_{\complex})\to \Omega_{\complex}$. Note that $\Omega_{\complex}\subseteq 
S^{G}_{\complex}(\gamma^*)$. Therefore, by part~\eqref{p1:local-diff}, the map
$\fKC$ is a local holomorphic diffeomorphism at every point in
$\fKC^{-1}(\Omega_{\complex})$. Moreover, the set $\fKC^{-1}(\Omega_{\complex})$ is compact.
Therefore, by~\cite[Proposition 2.19]{JML:03}, the map $\fKC\mid_{\fKC^{-1}(\Omega_{\complex})}$ is a covering map. Since the set $\Omega_{\complex}$ is simply
connected, by~\cite[Proposition A.28]{JML:03}, the mapping $\fKC\mid_{\fKC^{-1}(\Omega_{\complex})}$ is a holomorphic
diffeomorphism. Part~\eqref{p3:invertible_analytic} of the theorem then follows from this result. 
\end{proof}

\bibliographystyle{plainurl}
\bibliography{alias,Main,FB}

\end{document}